
\magnification=1200
\input amstex
\documentstyle{amsppt}
\def\Z{{\Bbb Z}}
\def\U{{\Cal U}}
\def\V{{\Cal V}}
\def\W{{\Cal W}}
\def\X{{\Cal X}}

\def\A{{\Cal A}}
\def\B{{\Cal B}}
\def\O{{\Cal O}}
\def\J{{\Cal J}}
\def\Nc{{\Cal N}}
\def\Zc{{\Cal Z}}
\def\Sc{{\Cal S}}

\def\Lc{{\Cal L}}
\def\M{{\Cal M}}
\def\D{{\Cal D}}
\def\coor#1#2{(#1_1,\allowmathbreak\dots,#1_{#2})}
\def\sh#1#2{\hbox{${\Cal #1}  #2 $}}

\def\Div{\operatorname{Div}}

\def\Hom{\operatorname{Hom}}
\def\Ker{\operatorname{Ker}}

\def\Id{\operatorname{Id}}
\def\Spec{\operatorname{Spec}}
\def\det{\operatorname{det}}
\def\medwedge{\textstyle\bigwedge}
\def\gr #1{\left\vert #1\right\vert}

\def\iso{\kern.35em{\raise3pt\hbox{$\sim$}\kern-1.1em\to}
         \kern.3em}
\def\sign{\operatorname{sign}}
\def\pd#1,#2{\displaystyle{\partial#1\over\partial#2}}
\def\volume{\left [ dz\otimes\pd ,\theta\right ]}
\def\rest #1,#2{{#1}_{\vert #2}}
\def\proof{\demo{Proof}}
\catcode`\@=11
\def\cub{\vrule height 1.4ex width 1.4ex depth -.1ex}
\def\qed{\ifhmode\unskip\nobreak\fi\ifmmode\ifinner\else\hskip5
\p@\fi\fi
 \hfill\hbox{\hskip5\p@\cub\hskip\p@}}
\catcode`\@=\active


\TagsOnRight
\topmatter
\title
THE VARIETY OF POSITIVE SUPERDIVISORS\\ OF A
SUPERCURVE (SUPERVORTICES)
\endtitle
\author J\.A\. Dom\'\i nguez
P\'erez, D\. Hern\'andez Ruip\'erez \& C. Sancho de Salas
\endauthor
\leftheadtext{J\.A\. Dom\'\i nguez, D\. Hern\'andez \& C. Sancho}
\rightheadtext{THE VARIETY OF POSITIVE SUPERDIVISORS\dots}
\affil
Departamento de Matem\'aticas Puras y Aplicadas, Universidad de
Salamanca \endaffil \address
Departamento de
Matem\'aticas Puras y Aplicadas, Universidad de Salamanca, Plaza de
la Merced 1-4, 37008 Salamanca, Spain. Phone: 3423294459. Fax:
3423294583 \endaddress
\email
sanz\@ relay.rediris.es (subject: To D.H.Ruiperez)
\endemail
\thanks
The first and the second authors acknowledge support received
under C\.I\.C\.Y\.T\. project PB-88-0379. The third author
was partially supported by D\.G\.I\.C\.Y\.T\. project
PS-88-0037\endthanks \subjclass 14A22, 14M30, 14D25, 14C05
\endsubjclass \abstract
The supersymmetric product of
a supercurve is constructed with the aid of a theorem of algebraic
invariants and the notion of positive relative superdivisor
(supervortex) is introduced. A supercurve of positive
superdivisors of degree 1 (supervortices of vortex number 1) of the
original supercurve is constructed as its supercurve of
conjugate fermions, as well as the supervariety of relative
positive superdivisors of degre $p$ (supervortices of vortex number
$p$.) A universal superdivisor is defined and it is proved
that every positive relative superdivisor can be obtained in a
unique way as a pull-back of the universal superdivisor. The case of
SUSY-curves is discussed.
\endabstract
\endtopmatter
\document
\heading
1. Introduction \endheading
Positive divisors of degree $p$ on an algebraic curve $X$ can be thought as
unordered sets of $p$ points of $X$, hence, as elements of the symmetric
$p$-fold product $S^pX$. The symmetric
$p$-fold product is the orbit space of the cartesian
$p$-fold product $X^p$ under the natural action of the symmetric group, and it
is thus endowed with a natural structure of algebraic variety. In this way,
positive divisors of degree $p$ are the points of an algebraic
variety $\Div^p(X)$, and this variety is of great importance in the
study of the geometry of curves, and it also has a growing interest
in Mathematical Physics.

{}From the geometrical side, one has, for instance, the role played
by the variety of positive divisors of degree $p$ in some classical
constructions
of the Jacobian variety of a complete smooth algebraic curve. The first
construction of the Jacobian variety, due to Jacobi and Abel, is of an analytic
nature and defines the Jacobian as a complex torus through the periods matrix.
The first algebraic construction is due to Weil \cite{31} who showed
that the algebraic structure and the group law of the Jacobian come
from the fact that it is birationally equivalent to the variety of
positive divisors of degree the genus of the curve. Another
procedure stemmed from Chow \cite{7}, who took advantage of the
fact that for $p$ high enough the Abel map (that maps a divisor of
degree $p$ into its linear equivalence class) is a projective
bundle, to endow the Jacobian with a structure of projective
algebraic group. But regardless the method used for constructing
the Jacobian, the structure of the variety of positive divisors of
degree $p$ and the diverse Abel morphisms from these varieties to
the Jacobian, turns out to be a key point in the theory of Jacobian
varieties (see, for instance \cite{17}, \cite{25}) and has proved
to be an important tool in the solution of the Schottcky problem
\cite{26}.

{}From a physical point of view, the variety of positive divisors of a complex
complete smooth curve $X$ (a compact Riemann surface) is the variety of
{\sl vortices\/} or solutions to the vortex equations (\cite{5},
\cite{10}.) For every holomorphic line bundle $L$ on
$X$ endowed with a hermitian metric, there is a Yang--Mills--Higgs
functional $YMH_\tau(\nabla,\phi)$ defined on gauge equivalence
classes of pairs $(\nabla,\phi)$ where $\nabla$ is an unitary
connection, by
$$
YMH_\tau(\nabla,\phi)=\int(\gr{F_\nabla}^2+
\gr{\nabla\phi}^2+\frac14\gr{\phi\otimes\phi\ast-\tau\Id}^2) d\mu
$$
where $F_\nabla$ is
the curvature, $\nabla\phi$ the covariant
derivative, and $\tau$ is a real parameter (see \cite{5}.)

Bradlow's theorem states that for large $\tau$, gauge
equivalence classes of solutions $(\nabla,\phi)$ to the vortex equation
$$
YMH_\tau(\nabla,\phi)=2\pi p\tau\,,
$$
where $p$ is the degree of $L$ with respect to the K\" ahler form,
correspond to divisors of degree $p$ on $X$. In this
correspondence, a solution $(\nabla,\phi)$ corresponds to the
divisor given by the set of centres of the vortices appearing with
multiplicity given by the multiplicity of the magnetic flux.

There is no similar theory for supersymmetric extensions of the vortex
equations
(supervortices,) and in fact only very little work on supervortices
or supersymmetric extensions of the Bogonolmy equations has been done (see
\cite{20}.) This paper will provide a first step in that direction, by
providing the right supervariety of positive superdivisors or supervortices on
a
supercurve.

This paper is organized as follows:

The {\sl supersymmetric product\/} $S^p\X$ for a supercurve $\X$ of dimension
$(1,1)$ is constructed in Part 2 as the orbit ringed space obtained through the
action of the symmetric group on the cartesian $p$-fold product of
$\X$. It is far from trivial that the resulting graded ringed space is a
supervariety of dimension $(p,p)$, a statement which is shown to be  equivalent
to an invariant theorem. It should be stressed that this
theorem is no longer true for supercurves of higher odd dimension, but our
result covers the most important cases such as SUSY-curves.

In Part 3 the notion of {\sl positive relative superdivisor\/} of degree $p$
for
a relative supercurve $\X\times\Sc\to\Sc$ is given. The classical definition
cannot be extended straightforwardly to supercurves if we wish that
superdivisors could be obtained as pull-backs of a suitable universal
superdivisor.

For ordinary algebraic curves, positive divisors of
degree 1 are just points. The novelty here is that
for an algebraic supercurve $\X$, positive relative superdivisors
of degree 1 (supervortices of vortex number 1) are {\sl are not\/}
points of $\X$ (see  \cite {Ma3},) but rather they are points of
another supercurve $\X^c$ with the same underlying ordinary
(bosonic) curve. Actually, if we think of $\X$ as a field of
fermions on a bosonic curve, the supercurve $\X^c$ is the
supercurve of conjugate fermions on the underlying bosonic curve.

This is proven in Part 4, that also contains the representability
theorem for positive relative superdivisors of degree $p$ on a
supercurve. The theorem means that positive relative superdivisors
of degree $p$ are the points of the supersymmetric $p$-fold product
$S^p \X^c$ of the supercurve $\X^c$ of conjugate fermions. This
property is stated in the spirit of Algebraic Geometry in terms of
the functor of the points: the precise statement is that the
functor of the positive relative superdivisors of degree $p$ of a
supercurve $\X$ of odd dimension 1, is the functor of the points of
the supersymmetric $p$-fold product $S^p \X^c$. This means that
every positive relative superdivisor of $\X\times\Sc\to\Sc$ can be
obtained in a unique way as the pull-back of a certain universal
positive superdivisor through a morphism $\Sc\to S^p \X^c$. We
obtain in that way what is the right structure of algebraic
superscheme the `space' of positive superdivisors of degree $p$ on
a supercurve can be endowed with.

The case of supersymmetric curves (SUSY-curves) is particularly important, both
by historical and geometrical reasons. We prove that for a supercurve $\X$, the
existence of a conformal stucture on $\X$ is equivalent to the existence of an
isomorphism between $\X$ and the supercurve $\X^c$ of conjugate fermions. In
other words, a supercurve $\X$ is a SUSY-curve, if and only if, $\X$ is
isomorphic with  $\X^c$. In this case the universal positive superdivisor of
degree 1 is Manin's superdiagonal (\cite{4}, \cite{23}) and we recover from
a clearer and more general viewpoint Manin's interpretation of the relationship
between points and positive divisors of degree 1 for SUSY-curves, and some
connected definitions (\cite{28}, \cite{29}.)

Summing up, the space of supervortices of vortex number $p$
(positive superdivisors of degree $p$) on a supercurve $\X$ of odd
dimension 1, is an algebraic supervariety of dimension $(p,p)$.
This algebraic supervariety is the supervariety $S^p \X^c$ of
`unordered families' of $p$ conjugate fermions. Moreover, only for
SUSY-curves supervortices of vortex number $p$ are `unordered
families' of $p$ points of $\X$.

This theory can be extended straightforwardly to SUSY-families
parametrized by a {\sl ordinary\/} algebraic scheme.

The results of this paper only in the case of SUSY-curves were
stated (without proofs) in \cite{8}.

\heading 2. Supersymmetric products\endheading
\subheading{1. Definitions}

A suitable reference for schemes theory is \cite{14};
the general theory of schemes in the supergeometry (superschemes)
can be found in \cite{22} and \cite{27}.

 Let $\X=(X,\A)$ be a graded ringed space, that is, a
pair consisting of a topological space $X$ endowed with a sheaf
$\A$ of $\Z_2$-graded algebras. Let us denote by $\J$ the ideal
$\A_1+\A_1^2$.
\definition{Definition 1} A superscheme of dimension $(m,n)$ over a
field $k$, is a graded ringed space $\X=(X,\A)$ where $\A$ is a
sheaf of graded $k$-algebras such that: \roster
 \item $(X,\O=\A/\J)$ is an $m$-dimensional scheme of finite type
over $k$. \item $\J/\J^2$ is a locally free $\O$-module
of rank $n$ and $\A$ is locally isomorphic to
$\bigwedge_{\O}(\J/\J^2)$. \endroster   \enddefinition
\definition{Definition 2}  A superscheme $\X=(X,\A)$ over a field
$k$ is said to be affine if the underlying scheme $(X,\O=\A/\J)$ is
an affine scheme, that is, if there is a homeomorphism
$$
X\iso\Spec(\Gamma(X,\O))
$$
and $\O$ is the sheaf on $X$ defined by localization on the basic
open subsets of the spectrum.\enddefinition
If $\X=(X,\A)$ is an affine superscheme, and $A=\Gamma(X,\A)$, we
shall simply write $\X=\Spec A$ for it.

 Let us consider the product $$\X^g=(X^g,\A^{\otimes g}) $$ where
$X^g$ denotes the cartesian product $X\times{\buildrel
g)\over\dots}\times X$, and $\A^{\otimes g}=\A\otimes{\buildrel
g)\over\dots}\otimes\A$.

The symmetric group $S_g$ acts on $\X^g$ by graded automorphisms of
superschemes according to the rule
$$
\aligned
\sigma\colon
X^g&\to X^g\\ (x_1,\dots,x_g)&\mapsto (x_{\sigma
(1)},\dots,x_{\sigma (g)})\\ \\
 \sigma^\ast\colon \A^{\otimes g}&\to
\sigma_\ast\A^{\otimes g}\\f_1\otimes\dots\otimes f_g
&\mapsto \prod \Sb i<j\\ \sigma(i)>\sigma
(j)\endSb (-1)^{\gr{f_i}\gr{f_j}} f_{\sigma(1)}\otimes\dots\otimes
f_{\sigma(g)}
\endaligned\tag 1
$$
where $\gr{\ }$
stands for the $\Z_2$-degree. This action reduces to the ordinary
action of $S_g$ on the scheme $(X^g,\O^{\otimes g})$. Then, we have
the orbit space $S^gX$, a natural projection $p\colon X^g\to S^gX$,
and an invariant sheaf $\O_g=\O^{S_g}$ on $S^gX$, whose sections  on
an open subset $V\subseteq S^gX$ are $$ \O_g(V)=\{f\in\O^{\otimes
g}(p^{-1}(V))\,\vert\,\sigma^\ast f=f {\hbox{\ for every\
}}\sigma\in S_g\}\,.$$

It is well-known that if $(X,\O)$ is a projective scheme, the
ringed space $(S^gX,\O_g)$ is a scheme, the {\sl symmetric
$p$-fold product\/} of $(X,\O)$ (\cite{30}, Prop.19.)

Let us consider the sheaf $\A_g=(\A^{\otimes g})^{S_g}$ of graded
invariants on $S^gX$ defined as above by letting $$
\A_g(V)=\{f\in\A^{\otimes g}(p^{-1}(V))\,\vert\,\sigma^\ast f=f
{\hbox{\ for every\ }}\sigma\in S_g\} $$ for every open subset
$V\subseteq S^gX$.

\subheading{2. The case of supercurves}
\definition{Definition 3}  A
supercurve is a superscheme $\X$ of dimension $(1,n)$ over a field
$k$. \enddefinition

 Let $\X$ be a  smooth proper supercurve, that is, a supercurve such
that $(X,\O)$ is proper and smooth.

\proclaim{Theorem 1}  Let $\X=(X,\A)$ be a smooth proper supercurve
of odd dimension $n>0$. The graded ringed space $S^g\X$ is a
superscheme if and only if $n=1$, that is, if and only if $\X$ is a
superscheme of dimension $(1,1)$. In that case, $S^g\X=(S^gX,\A_g)$
is a superscheme of dimension $(g,g)$ that will be called the
supersymmetric $g$-fold product of $\X$.
\endproclaim
\proof Let us notice that $(X,\O)$ is projective
 (it has very ample sheaves,) so that the ringed space $(S^gX,\O_g)$
is a scheme as we mentioned above (in fact, it is smooth, which is
no longer true for higher dimensional $X$.)

As there is a natural
projection $\A_g\to\O_g$, we have only to ascertain if $\A_g$ is
locally the exterior algebra of a locally free $\O_g$-module. We
can thus assume $\A=\bigwedge_{\O}(\Nc)$, $\Nc$ being a free rank
$n$ $\O$-module.

Let us write $\Nc_i=\O\otimes\dots\otimes{\buildrel{\phantom{i)}
\downarrow i)}\over\Nc}\otimes\dots\otimes\O$ and
$\M=\Nc_1\oplus\dots\oplus\Nc_g$. Now, if  $\bar\O=\O^{\otimes g}$
and $\bar\A=\A^{\otimes g}$, we have $$ \bar\A={\textstyle
\bigwedge}_{\O}(\Nc)\otimes_{\O}
{\buildrel{g)}\over\dots}\otimes_{\O}{\textstyle
\bigwedge}_{\O}(\Nc)\iso {\textstyle \bigwedge}_{\bar\O}(\M)\,. $$

The symmetric group $S_g$ acts on $\M$ by $$
\sigma(n_1+\dots+n_g)=n_{\sigma(1)}+\dots+n_{\sigma(g)} $$ and this
action provides an action $\sigma\colon\bar\A\to\bar\A$ on the
exterior algebra $\bar\A=\bigwedge_{\bar\O}(\M)$, given by $$
\sigma(m_1\wedge\dots\wedge
m_p)=\sigma(m_1)\wedge\dots\wedge\sigma(m_p) $$

This action of $S_g$ on $\bar\A$ is actually equal to the one
defined in (1), because both coincide on
$\bigwedge_{\bar\O}^1(\M)=\M$ and are morphisms of graded algebras.

If we denote by $\M^{S_g}$ the $\O_g$-module consisting of the
invariant sections of $\M$, the proof of Theorem 1
will be
thus completed with the following \proclaim{Lemma 1}  The natural
morphism of sheaves of graded $\O_g$-algebras over $S^gX$, $$
\phi\colon\medwedge_{\O_g}(\M^{S_g})\to
(\medwedge_{\bar\O}(\M))^{S_g}=\A_g\,, $$ is an isomorphism if and
only if $n=1$.\endproclaim \proof The proof is a
computation of invariants in the exterior algebra of a free module
over a {\sl commutative\ } ring, which allows us to use standard
methods of Commutative Algebra (all the results that we shall use
can be found, for instance, in \cite{1}.)

Let us start with the case $n=1$.

a) We can assume that $X=\Spec\O$, where $\O$ is a semilocal ring
with $g$ maximal ideals $\frak p_1,\dots,\frak p_g$, and then, that
$\Nc=\O\cdot e$, $\Nc_i=\bar\O\cdot e_i$ (where
$e_i=1\otimes\dots\otimes{\buildrel{\phantom{i)} \downarrow
i)}\over e}\otimes\dots\otimes1$) and $\M=\bar\O\cdot
e_1\oplus\dots\oplus\bar\O\cdot e_g$.

Let us notice, first, that $\phi$ is an isomorphism if and only if
it is an isomorphism when localized at every maximal ideal
$\frak{p}$ of $\O_g$. On the other hand, $\frak p$ corresponds to a
divisor $D=x_1+\dots+x_g$ and  the fibre of $p\colon X^g\to S^gX$
over this point consists of the family $\coor xg$ (some of the
points can be equal) together with its permutations. It follows
that we are reduced to consider only the localization of $\O$ at
these particular points $\coor xg$.

b) We can assume that $\O=k[t]$ and $\Nc=k[t]\cdot e$.

Since the completion morphism $(\O_g)_{\frak{p}}\hookrightarrow
\widehat{(\O_g)}_{\frak{p}}$ is a faithfully flat morphism, we are
reduced to show that $\phi$ is an isomorphism after completing
$\O_g$ at every maximal ideal.

Let $t\in\O$ be an element that takes different values $\coor
\lambda g$ at the points $\coor xg$ and such that $t-\lambda_i$ is
a parameter at $x_i$ (that is, it generates the maximal ideal of
the local ring $\O_{{\frak p}_i}$.) Then, $k[t]$ is a subring of
$\O$ and  moreover, given two different maximal ideals $\frak
p_i\ne\frak p_j$, the maximal ideals $\bar{\frak p}_i=\frak p_i\cap
k[t]$ and $\bar{\frak p}_j=\frak p_j\cap k[t]$ of $k[t]$ are also
different.

Let $I=\frak p_1\cap\dots\cap\frak p_g$, $\bar I=\bar{\frak
p}_1\cap\dots\cap\bar{\frak p}_g$ be the intersection ideals and
$\O\hookrightarrow\widehat\O$, $k[t]\hookrightarrow\widehat{k[t]}$
the faithfully flat morphisms of completion with respect to the
ideals $I$, $\bar I$, respectively. Then we have $$
\widehat{k[t]}\iso\prod_{i=1}^g\widehat{k[t]}_{{\frak p}_i}\iso
\prod_{i=1}^g\widehat{\O}_{{\frak p}_i}\iso\widehat\O\,, $$ so
that, if the theorem is true for $k[t]$ and the module $k[t]\cdot
e$, it is also true for $\widehat{k[t]}=\widehat\O$ and
$\widehat{k[t]}\cdot e=\widehat\O\cdot e$ and then it will be true
for $\O$ and $\Nc=\O\cdot e$ by faithfull flatness.

c) The case $\O=k[t]$ and $\Nc=k[t]\cdot e$.

Now, for every $0\leq p\leq g$ we have $$
\medwedge^p(\M)=\bigoplus_{i_1<\dots<i_p}
\Nc_{i_1}\wedge\dots\wedge\Nc_{i_p}
 $$ and $S_g$ acts transitively by permutation of summands.

In fact, $\Nc_{i_1}\wedge\dots\wedge\Nc_{i_p}= \sigma_{i_1\dots
i_p}(\Nc_1\wedge\dots\wedge\Nc_p)$, $\sigma_{i_1\dots i_p}\in S_g$
being any permutation of type $\sigma_{i_1 \dots i_p}= \pmatrix
1&\dots&p&\dots\\ {i_1}&\dots&{i_p}&\dots \endpmatrix$.

Then, an invariant element
$m=\sum_{i_1<\dots<i_p}n_{i_1}\wedge\dots\wedge n_{i_p}$ is
characterized by $n_1\wedge\dots\wedge n_p$ and we have an
isomorphism $$\aligned (\Nc_1\wedge\dots\wedge\Nc_p)^{S_p\times
S_{g-p}}&\iso(\medwedge^p(\M))^{S_g}\\ n_1\wedge\dots\wedge
n_p&\mapsto\sum_{\sigma\in S_g/(S_p\times S_{g-p})} \sigma
(n_1\wedge\dots\wedge n_p) \endaligned
 $$ where $S_p\times S_{g-p}$ denotes the subgroup of $S_g$
consisting of the permutations leaving invariant the subset
$\{1,\dots,p\}$. In particular, $(\Nc_1)^{S_{g-1}}\iso\M^{S_g}$.

Since $\Nc_1\wedge\dots\wedge\Nc_p=\bar\O\cdot e_1\wedge\dots\wedge
e_p$,
 we have $$ (\Nc_1\wedge\dots\wedge\Nc_p)^{S_p\times
S_{g-p}}\iso\bar\O^{-S_p\times S_{g-p}}
 $$ where $\bar\O^{-S_p\times S_{g-p}}$ stand for the subset of
those $f\in\bar\O$ such that $(\sigma\times\mu)(f)=
\sign(\sigma)\cdot f$ for every $(\sigma\times\mu)\in S_p\times
S_{g-p}$. Taking $p=1$, we obtain $$
\M^{S_g}\iso\Nc_1^{S_{g-1}}\iso\bar\O^{S_{g-1}}\,, $$ and the
original morphism $$\phi\colon\medwedge_{\O_g}^p (\M^{S_g})\to
(\medwedge_{\bar\O}^p(\M))^{S_g}=\A_g
 $$ is now the morphism $$
\bar\phi_p\colon\medwedge^p\bar\O^{S_{g-1}} \to\bar\O^{-S_p\times
S_{g-p}}
 $$ described by $$ \bar\phi_p(f_1\wedge\dots\wedge
f_p)=\sum_{\mu\in
S_p}\sign(\mu)\sigma_{\mu(1)}(f_1)\dots\sigma_{\mu(p)}(f_p)
$$
where $\sigma_i$ is the transposition of 1 and $i$.

Let us proof that $\bar\phi_p$ is an isomorphism: There is a
commutative diagram
$$
\CD
\bar\O^{S_{g-1}}\otimes_{\O_g}\dots\otimes_{\O_g}\bar\O^{S_{g-1}}
@>T>>\bar\O^{1\times S_{g-p}}\\
@V H VV @VV H' V\\
\medwedge_{\O_g}^p\bar\O^{S_{g-1}}
@>{\bar\phi_p}>>\bar\O^{-S_p\times S_{g-p}}
\endCD
$$
where
$T(f_1\otimes\dots\otimes f_p)=\sigma_1(f_1)\dots\sigma_p(f_p)$,
$H(f_1\otimes\dots\otimes f_p)=f_1\wedge\dots\wedge f_p$ and
$H'(f)=\sum_{\mu\in S_p}\sign(\mu)(\mu\times 1)(f)$.

As $\O=k[t]$, $\bar\O=k[t_1,\dots,t_g]$, and if we denote
$$\align(s_1,\dots,s_g)&={\hbox{\rm symmetric functions of\
}}(t_1,\dots,t_g)\\ (\bar s_1,\dots,\bar s_{g-1})&={\hbox{\rm
symmetric functions of\ }}(t_2,\dots,t_g)\\
(s'_1,\dots,s'_{g-p})&={\hbox{\rm symmetric functions of\ }}
(t_{p+1},\dots,t_g)\\ \intertext{we have} \O_g&=k[s_1,\dots,s_g]\cr
\bar\O^{S_{g-1}}&=k[t_1,\bar s_1,\dots,\bar
s_{g-1}]\iso\O_g[t_1]\cr \bar\O^{1\times
S_{g-p}}&=k[t_1,\dots,t_p,s'_1,\dots,s'_{g-p}]\iso
\O_g[t_1,\dots,t_p]\,. \endalign $$ If follows that
$\bar\O^{-S_p\times S_{g-p}}$ can be identified with the
$p^{\hbox{{\sevenrm th}}}$ skew-symmetric tensors of the
$\O_g$-module $\bar\O^{S_{g-1}}=\O_g[t_1]$ and the previous diagram
reads
$$
\CD
\O_g[t_1]\otimes_{\O_g}{\buildrel
p)\over\dots}\otimes_{\O_g}\O_g[t_1]@>\sim>>
\O_g[t_1,\dots,t_p]\\
@V H VV @VV H' V\\
\medwedge_{\O_g}^p\O_g[t_1]@>{\bar\phi_p}>>
\medwedge_{\O_g}^p\O_g[t_1] \endCD $$ where now $H'$ is the
skew-symmetrization operator, finishing the proof of the if part.

To complete the proof, we have to show that if $n>1$, $\bar\phi_p$
is not an isomorphism. Let us write $\Nc=\bigoplus_{j=1}^{n}\Nc^j$
with $\Nc^j$ of rank 1, and
$\M^j=\bigoplus_{i=1}^{g}(\O\otimes\dots\otimes
{\buildrel\phantom{i)}\downarrow
i)\over{\Nc^j}}\otimes\dots\otimes\O)$ so that
$\M=\bigoplus_{j=1}^n\M^j$.

Then $$\aligned (\medwedge_{\bar\O}(\M))^{S_g}&=
\bigoplus_{p_1+\dots+p_s=p}
(\medwedge_{\bar\O}^{p_1}\M^1\otimes\dots\otimes
\medwedge_{\bar\O}^{p_s}\M^s)^{S_g}\\ \medwedge_{\O_g}(\M^{S_g})&=
\bigoplus_{p_1+\dots+p_s=p}
\medwedge_{\O_g}^{p_1}(\M^1)^{S_g}\otimes\dots\otimes
\medwedge_{\O_g}^{p_s}(\M^s)^{S_g}\,. \endaligned $$

By the case $n=1$, we have $$ \medwedge_{\O_g}^{ p_i
}(\M^i)^{S_g}\iso({ \textstyle\bigwedge}_{\O_g}^{ p_i
}(\M^i))^{S_g} $$ and then $$ \medwedge_{\O_g}(\M^{S_g})=
\bigoplus_{p_1+\dots+p_s=p}
(\medwedge_{\O_g}^{p_1}\M^1)^{S_g}\otimes\dots\otimes
(\medwedge_{\O_g}^{p_s}\M^s)^{S_g}\,. $$

But there are invariant elements in the tensor product
$(\medwedge_{\bar\O}^{p_1}\M^1 \otimes\dots\otimes
\medwedge_{\bar\O}^{p_s}\M^s)^{S_g}$ which cannot be written as
tensor products of invariant elements. This means that the morphism
$\phi_p$ is not an isomorphism in this case. \qed\enddemo\enddemo
\proclaim{Corollary 1} If $(z,\theta)$ are graded local coordinates
on a supercurve $\X$ of dimension $(1,1)$, a system of graded local
coordinates for $S^g\X$ is given by
$(s_1,\dots,s_g,\varsigma_1,\dots,\varsigma_g)$
$(s_1,\dots,s_g)$ are the (even) symmetric functions of
$(z_1,\dots,z_g)$ and $(\varsigma_1,\dots,\varsigma_g)$ are the odd
symmetric functions defined by
$\varsigma_h=\sum_{i=1}^g\sigma_i(\theta_1\bar
s_{h-1})$.\endproclaim
 \heading 3. Positive superdivisors\endheading
 From this point, calligraphic
types are reserved to graded ringed spaces and the structure ring
sheaf of any ringed space will be denoted by $\O$ with the name of
the ringed space as a subscript. For instance, $\X=(X,\O_\X)$ or
simply $\X$ will mean a graded ringed space, whereas $(X,\O_X)$ or
$X$ will represent the underlying ordinary ringed space.

\subheading{1. The universal divisor for an algebraic curve}

 This section is devoted to summarize the theory of the variety of
positive divisors and the universal divisor for a (ordinary) smooth
proper algebraic curve $X$, and to show that the universal
property still holds when the space of parameters is a
superscheme. Suitable
references are \cite{13} or \cite{16}.

In that case, positive divisors of degree $g$ are unordered
families of $g$ points, and they are then parametrized by the space
of such families, that is, by the symmetric product $S^gX$. This
can be made precise through the notion of relative divisor.

If $S$ is another scheme, positive relative divisors of $X\times
S\to S$ of degree $g$ are subschemes $Z\to X$ such that $\O_Z$ is a
locally free $O_S$-module of rank $g$. There is a nice positive
relative divisor $Z^u$ of degree $g$ of $X\times S^gX\to S^gX$,
whose fibre on a point $(x_1,\dots,x_g)\in S^gX$ is the divisor
$x_1+\dots+x_g$ of $X$ defined by it. $Z^u$ is called the {\sl
universal divisor\/} because the map $$\aligned
\Hom(S,\,S^gX)&\to\Div_S^g(X\times S)\\
\phi&\mapsto(1\times\phi)^{-1}(Z^u) \endaligned\, $$ where
$\Div_S^g(X\times S)$ denotes the set of positive relative divisors
of degree $g$, is {\sl one to one\/}. This means that each positive
divisor can be obtained as a pull-back of the universal divisor;
this statement is known as {\sl representability theorem\/} for
the symmetric product.

But it turns out that the above theory is still true when a
superscheme is allowed as the space of parameters, once the
corresponding notion of positive relative divisor has been
established.
\definition{Definition 4} Let $X$ be an ordinary smooth
curve and $(\Sc,\O_{\Sc})$ a superscheme. A positive relative
divisor of degree $g$ of $X\times\Sc\to\Sc$ is a closed
sub-superscheme $\Zc$ of $X\times\Sc$ of codimension $(1,0)$
defined by a homogeneous ideal $J$ of $\O_{X\times\Sc}$ such that
$\O_{X\times\Sc}/J$ is a locally free $\O_{\Sc}$-module of rank
$(g,0)$.\enddefinition The ideal $J$ of a positive
relative divisor of degree $g$ is then locally generated by an
element of type $$
f=z^g-a_1z^{g-1}+\dots+(-1)^ga_g\,,\tag 2
$$ where the
$a_i$'s are even elements in $\O_{\Sc}$, and $\O_{X\times\Sc}/J$ is
a free $\O_{\Sc}$-module with basis $(1,z,\dots,z^{g-1})$.

 The representability theorem now reads \proclaim{Theorem 2}   Let
$X$ be a smooth proper curve over a field $k$ and $Z^u$ the
universal divisor. The map $$ \eqalign{
\Hom(\Sc,\,S^gX)&\to\Div_{\Sc}^g(X\times \Sc)\cr
\phi&\mapsto(1\times\phi)^{-1}(Z^u)\,,}\tag 3 $$ where
$\Div_{\Sc}^g(X\times \Sc)$ denotes the set of positive relative
divisors of degree $g$, is one to one for every superscheme $\Sc$.
\qed
\endproclaim
 Proof of the representability theorem for ordinary schemes applies
with only minor changes to this case.

There are two key points for the proof of this theorem. The first
one is the construction of the universal divisor, which can be done
as follows: If $\pi_i\colon X^g\to X$ is the ith projection and
$\Delta_i$ is the positive relative divisor of $X\times X^g\to X^g$
obtained by pull-back of the diagonal $\Delta\subset X\times X$
throughout $1\times\pi_i\colon X\times X^g\to X\times X$ we can
prove there is a unique positive relative divisor $Z^u$ of $X\times
S^gX\to S^gX$ such that $$ (1\times p)^{-1}
Z^u=\Delta_1+\dots+\Delta_g $$ where $p\colon X^g\to S^gX$ is the
natural projection. This divisor $Z^u$ is the universal divisor.

The second key point is the so-called `determinant morphism'
$\Sc\to S^g \Zc$, $\Zc$ being a positive relative divisor of degree $g$
because its composition with $S^g \Zc\to S^gX$ provides the inverse
mapping of (3) (See \cite{Iv}.) The determinant
morphism for the locally free $\O_{\Sc}$-module of rank $(g,0)$
$\O_\Zc$ is defined as follows: Each element $b$ in the invariant
sheaf $(\O_\Zc)_g=(\O_\Zc^{\otimes g})^{S_g}$ acts on the
$\O_{\Sc}$-module $\bigwedge_{\O_\Sc}\O_\Zc$ of rank $(1,0)$ as the
multiplication by a well-determined element $\det (b)$ in
$\O_{\Sc}$. This gives rise to a morphism of sheaves
$(\O_\Zc)_g\to\O_{\Sc}$, and to a morphism of schemes $\Sc\to S^g
\Zc$. The determinant morphism provides the inverse mapping of
(3) because if $b$ is an even element in $\O_\Zc$,
$b_i=1\otimes\dots\otimes{\buildrel{\phantom{i)} \downarrow
i)}\over b}\otimes\dots\otimes1\in \O_\Zc^{\otimes g}$, and we
denote by $s_i(b)$ the symmetric functions of $b_1,\dots,b_g$, we
have that $$ a_i=\det(s_i(b))\quad (i=1,\dots,g)\,, $$ where
$z^g-a_1z^{g-1}+\dots+(-1)^ga_g$ is the characteristic polynomial
of $b$ acting on $\O_\Zc$ by multiplication (compare with
(2).)

\subheading{2. Positive superdivisors on supercurves}

The above discussion is based on a trivial but important point:
positive divisors are families of points. Even in the relative
case, positive relative divisors of degree 1 are `$S$-points,' that
is, sections of $X\times S\to S$, and by this reason, positive
divisors of degree $g$ are parametrized by the symmetric product
$S^gX$ and the universal divisor.

For a supercurve $\X=(X,\O_\X)$, a similar notion could be done, by
defining positive relative superdivisors of $\X\times\Sc\to\Sc$
($\Sc$ being an arbitrary superscheme,) as closed sub-superschemes
of $\X\times\Sc$ of codimension $(1,0)$ flat over
the base superscheme.

This definition has two drawbacks. The first one is that
`$\Sc$-points' are not superdivisors in that sense because they
have codimension $(1,1)$ and not codimension $(1,0)$ as
superdivisors does (\cite{23}.) The second one is that we cannot
ensure that they are pull-backs of a suitable universal
superdivisor.

We have thus modified the notion of positive relative superdivisors
in order to fulfill the second requirement as follows:

Let ($\X,\O_\X)$ be a  smooth supercurve, and $(\Sc,\O_{\Sc})$ a
superscheme.
 \definition{Definition 5}  A positive relative superdivisor
of degree $g$ of $\X\times\Sc\to\Sc$ is a closed sub-superscheme
$\Zc$ of $\X\times\Sc$ of codimension $(1,0)$ whose reduction
$\hat\Zc=\Zc\times_{\X}X$ is a positive relative divisor of degree
$g$ of $X\times\Sc\to\Sc$ (see Definition 4.) \enddefinition

Even with our definition, `$\Sc$-points' are not superdivisors, but
as we shall see afterwards, there is a close relationship between
them, at least for SUSY-curves.

 Positive relative superdivisors can be described locally in a
rather precise way in the case of a smooth supercurve of dimension
$(1,1)$.

In this case, the natural morphism $\O_{\X}\to\O_X$ induces an
isomorphism $(\O_{\X})_0\iso\allowmathbreak\O_X$, so that $\O_{\X}$
is a module over $\O_X$, there exists a canonical projection $\X\to
X$ and $\O_{\X}$ is in a natural way an exterior algebra
$\O_{\X}=\bigwedge_{\O_{\X}} \Lc$, where $\Lc=(\O_{\X})_1$ is a
line bundle over the ordinary curve $X$.

\proclaim{Lemma 2} Let $\X$ be a smooth supercurve of dimension
$(1,1)$. A closed sub-superscheme $\Zc$ of $\X\times\Sc$ of
codimension $(1,0)$ defined by a homogeneous ideal $J$ of
$\O_{\X\times\Sc}$ is a positive relative superdivisor of degree
$g$ if and only if the following conditions hold:
\roster
\item
$\O_\Zc=\O_{\X\times\Sc}/J$ is a locally free $\O_\Sc$-module of
dimension $(g,g)$.
\item If $(z,\theta)$ is a system of graded local
coordinates, $J$ can be locally generated by an
element of type $$ f=z^g-(a_1+\theta
b_1)z^{g-1}+\dots+(-1)^g(a_g+\theta b_g) $$ where the $a_i$'s are
even and the $b_j$'s are odd elements in $\O_\Sc$.
\endroster
\endproclaim \proof Let $\Zc$ be a positive relative superdivisor of
degree $g$ defined by a homogeneous ideal $J$ of
$\O_{\X\times\Sc}$ and let us consider a system of relative local
coordinates $(z,\theta)$. Then, the reduction
$\hat\Zc=\Zc\times_{\X}X$ is a positive relative divisor of degree
$g$ of $X\times\Sc\to\Sc$ defined by the image $\hat J$ of $J$ by
the morphism $\pi\colon\O_{\X\times\Sc}\to\O_{X\times\Sc}$, so
that an element $f\in J$ generates $J$ if and only if $\hat J$ is
generated by $\hat f=\pi(f)$. Since $\hat J$ defines a positive
relative divisor of degree $g$ of $X\times\Sc\to\Sc$,
 then $\hat J$ has a generator of
type $\hat f=z^g-a_1z^{g-1}+\dots+(-1)^ga_g$ where the $a_i$'s are
even elements in $\O_{\Sc}$ (see equation (2),) and
$\O_{\hat\Zc}=\O_{X\times\Sc}/\hat J$ is a free $\O_{\Sc}$-module
with basis $(1,z,\dots,z^g)$. This means that
$\O_{\hat\Zc}\iso\O_{\Sc}[z]/(\hat f)$. It follows that there is a
generator of $J$ of the form $f=\hat f+\theta\cdot d$ and that
$d\equiv q(z)\pmod{\hat J}$ for certain polynomial $q(z)$ of degree
less than $g$. In consequence, the element $\hat f+\theta q(z)$
generates $J$ and is of the predicted type. An easy computation now
shows that $\O_{\Zc}$ is a rank $(g,g)$ free $\O_{\Sc}$-module with
basis $(1,z,\dots,z^{g-1},\theta,\theta z,\dots,\theta z^{g-1})$.

The converse is straightforward. \qed\enddemo

\subheading{3. The functor of positive superdivisors on a
supercurve}

Let $(\X,\Zc)$ be a supercurve. For every superscheme $\Sc$ let
us denote by $\Div_\Sc^g(\X\times \Sc)$ the set of positive
relative superdivisors of degree $g$ of $\X\times\Sc\to\Sc$. If
$\varphi\colon\Sc'\to\Sc$ is a morphism of superschemes, and $\Zc$
is a positive relative divisor of degree $g$ of
$\X\times\Sc\to\Sc$, $(1\times\varphi)^{-1}\Zc$ is a positive
relative divisor of degree $g$ of $\X\times\Sc'\to\Sc'$. In
categorial language this essentially means that $$
\Sc\to\Div_\Sc^g(\X\times \Sc) $$ is a functor.

We whish to show that when $\X$ has dimension $(1,1)$, the above
functor is {\sl representable\/} in a similar sense to that of
Theorem 2. A proof is given in the next
section.
 \heading 4. The representability theorem for positive
superdivisors\\on a supercurve of dimension
$(1,1)$\endheading

In what follows, we consider only supercurves $\X=(X,\O_{\X})$
which are smooth, proper and of dimension $(1,1)$. This last
condition means that the structure sheaf $\O_{\X}$
is canonically isomorphic with $\O_X\oplus\Lc$ for certain
line bundle $\Lc$ on the ordinary underlying curve $X$.

\subheading{1. The supercurve of positive
superdivisors of degree 1}

Let $\Sc=(\Spec B,\B)$ an affine superscheme and
$\Zc=(Z,\O_{\Zc})\hookrightarrow \X\times \Sc\to \Sc$ a relative
superdivisor of degree 1. The structure sheaf $\O_{\Zc}$ is a
quotient of the structure sheaf $(\O_X\oplus\Lc)\otimes_k
\B$ of $\X\times\Spec B$. We also have that $\O_{\Zc}\iso \B\oplus
\overline {\Lc}$ where $\overline {\Lc}$ is the image of
$\Lc\otimes_k\B$ in $\O{_\Zc}$, since
$\O_{\hat\Zc}\iso\B$, because $\Zc$ is a superdivisor of
degree 1. Moreover,  $\overline {\Lc}=\Lc\otimes_{\O_X}\B$,
where $\B$ is an $\O_X$-algebra trough the natural morphism
$f\colon \O_X\to\O_{\hat\Zc}\iso\B$, so that it is a locally
free rank 1 $\B$-module.

It is now clear that the superdivisor $\Zc$ is characterized
by the morphism
$f\colon \O_X\to\B$ together with a morphism
$\tilde f\colon \O_X\to \B\oplus \overline {\Lc}$ extending
$f$. That is, $\Zc$ is defined by a morphism $f\colon \O_X\to\B$
and a derivation $\Delta\colon\O_X\to\overline
{\Lc}_0=\Lc\otimes_{\O_X}\B_1$. But $\Delta$ can be understood as
an element $f_{\Delta}\in \Hom_{\O_X}(\kappa_X, \overline
{\Lc}_0)\iso\Hom_{\O_X}(\kappa_X
\otimes_{\O_X}\Lc^{-1},\B_1)$ (where $\kappa_X$ is the
canonical sheaf of $X$,) so that the couple
$(f,\Delta)$ is equivalent to a graded ring morphism $g\colon
\O_X\oplus (\kappa_X\otimes_{\O_X}\Lc^{-1})\to
\B$.

The above discusion remains true for arbitrary (non affine)
superschemes $\Sc$. This means that the supercurve
$\X^c=\Spec (\O_X\oplus\Lc^c)$,
where $\Lc^c=\kappa_X\otimes_{\O_X}\Lc^{-1}$, will
represent the functor of superdivisors of degree 1 of
the supercurve $\Spec (\O_X\oplus\Lc)$. The
universal divisor, ${\Zc}^u_1\hookrightarrow\X\times\X^c$, will be
the divisor corresponding to the identity morphism
$\Id\colon \Sc=\X^c\to\X^c$. One can
compute this superdivisor as above and obtain that it
is the closed subsuperscheme whose ideal sheaf is the
kernel of the graded ring morphism:
$$
\bar\partial\colon
(\O_X\oplus\Lc)\otimes_k(\O_X\oplus\Lc^c)=
\medwedge_{\O_X\otimes_k\O_X}[(\Lc\otimes_k\O_X)\oplus
(\O_X\otimes_k \Lc^c)]\to\medwedge_{\O_X}(\Lc\oplus \Lc^c)
$$
given by $a\otimes b\mapsto a\cdot
b\oplus b\cdot d(a)$ on $\O_X\otimes_k\O_X$ (taking into
account that $b\cdot d(a)$ is a local section of
$\kappa_X\iso\Lc\otimes_{\O_X}\Lc^c$)
and as the natural morphisms on the remaining components.
Moreover, if  $U\subset X$ is an affine open subset and
$z\in \O_X(U)$ is a local parameter, and if
$\Lc$ is trivial on $U$, $\rest\Lc,U\iso\theta
\cdot\rest{\O_X},U$, then $\Lc^c$ is trivial on
$U$ generated by
${\theta}^c=\omega _{\theta}\cdot dz$,
$\omega_{\theta}\in\Gamma (U,\Lc^{-1})$ being the dual
basis of $\theta$. If $\U=\Spec
(\O_X\oplus\Lc)\subset \X$ and  ${\U}^c=\Spec
(\O_X\oplus\Lc^c)\subset \X^c$, the
restriction of the universal superdivisor ${\Zc}^u_1$ to
$\U\times{\U}^c$ is given by the local equation:
$$
z_1-z_2-\theta
\otimes\theta^c=0\tag 4
$$
where $z_1=z\otimes 1$ and $z_2=1\otimes z$.

The above discussion can be summarized as follows:

 Let $\X=(X,\O_{\X}=\O_X\oplus\Lc)$ be a smooth
proper supercurve of dimension $(1,1)$.
\definition{Definition 6}  The supercurve of positive divisors
of degree 1 on $\X$ is the supercurve of dimension (1,1)
defined as  $\X^c=(X,\O_X\oplus\Lc^c)$
where $\Lc^c=\kappa_X\otimes_{\O_X}\Lc^{-1}$. This
supercurve is also called the supercurve of conjugate fermions on
$\X$. \enddefinition \definition{Definition 7}  The universal
positive superdivisor of degree 1 is the relative superdivisor
${\Zc}^u_1$ of $\X\times\X^c\to\X^c$ defined by the ideal sheaf
$\Ker\bar\partial$ earlier considered. If $(z,\theta)$ are graded
local coordinates for $\X$, the corresponding local equation of
${\Zc}^u_1$ is $z_1-z_2-\theta
\otimes\theta^c=0$  where $z_1=z\otimes 1$, $z_2=1\otimes z$ and
${\theta}^c=\omega _{\theta}\cdot dz$.
\enddefinition \proclaim{Theorem 3} The morphism of functors:
$$ \aligned \Theta\colon
\Hom(\Sc,\X^c)&\to\Div_\Sc^1(\X\times \Sc)\\
\varphi&\mapsto (1\times\varphi)^{-1}(\Zc^u_1)\,. \endaligned
$$
is a functorial isomorphism.
\endproclaim
By this representability theorem, the supercurve
$\X^c$ of conjugate fermions para\-me\-trizes positive
superdivisors of degree 1 on the original supercurve $\X$. That
means that positive superdivisors of degree 1 on $\X$ are not points
of $\X$ as it happens in the ordinary case, but rather points
of another supercurve $\X^c$ with the same
underlying ordinary curve $X$.

\subheading{2. Positive superdivisors of degree 1 on a
SUSY-curve}

This section will explore the relationship between points and
positive superdivisors of degree 1 for a SUSY-curve
(Supersymmetric curve.)
This
relationship was firstly described by Manin (see \cite{23},) but
it can be enlightened by means of the supercurve of
positive superdivisors of degree 1 defined above. Let us start
by recalling some definitions and elementary properties of
SUSY-curves. More details can be found in Manin
\cite{21}, \cite{22}, \cite{23}, \cite{24}, Batchelor and
Bryant \cite{3}, Falqui and Reina \cite{9}, Giddings and Nelson
\cite{11}, \cite{12}, Bartocci, Bruzzo and Hern\'andez
Ruip\'erez \cite{2}, Bruzzo and Dom\'{\i}nguez P\'erez
\cite{6}, or LeBrun, Rothstein, Yat-Sun Poon and Wells
\cite{18}, \cite{19}.

Let $\Sc=(S,\O_\Sc)$ be a superscheme.
\definition{Definition 8}  A supersymmetric curve or SUSY-curve over
$\Sc$, is a proper smooth morphism
$\X=(X,\O_\X)\to\Sc$ of superschemes of relative dimension $(1,1)$
endowed with a locally free submodule $\D$ of rank $(0,1)$ of the
relative tangent sheaf $T_{\X/\Sc}=\sh Der_{\O_\Sc}(\O_\X)$ such
that the composition map $$\D\otimes_{\O_\X}\D\buildrel
[\;,\;]\over\longrightarrow\sh Der_{\O_\Sc}(\O_\X)\to \sh
Der_{\O_\Sc}(\O_\X)/\D$$ is an isomorphism of $\O_\X$-modules
(see, for instance, \cite{19}.) \enddefinition

If $\X=(X,\O_\X,\D)$ is a SUSY-curve, $X$ can be covered by affine
open subsets $U\subseteq X$ with local relative coordinates
$(z,\theta)$ such that $\D$ is locally generated by
$D=\pd,\theta+\theta\pd,z$. These coordinates are called {\sl
conformal\/}.

There is a natural isomorphism $\D^\ast\iso Ber_{\O_\Sc}(\O_\X)$
and a `Berezinian differential' $$
\partial\colon\Omega^1_{\X/\Sc}\to \D^\ast\simeq
Ber_{\O_\Sc}(\O_\X)$$ which is nothing but the natural projection
induced by the immersion $\D\to T_{\X/\Sc}$. In conformal
coordinates $\partial$ is described by $\partial(df)=\volume\cdot
D(f)$, where $\volume$ denotes the local basis of
$Ber_{\O_\Sc}\O_\X$ determined by $(z,\theta)$ (see \cite{15},
\cite{23}.)

If $(X,\O_{\X},\D)$ is a (single) SUSY-curve, that is, a SUSY-curve
over a point, we have that $\O_{\X}=\bigwedge_{\O_{\X}}(\Lc)$, and
there are isomorphisms $\D\otimes_{\O_{\X}}\O_X\iso\Lc^{-1}$ and
$$
\Lc\otimes_{\O_X}\Lc\iso\kappa_X\,.
$$
This isomorphism if often called a spin structure on $X$.
Conversely, a spin structure induces a conformal structure, so
that a conformal structure on a proper smooth
supercurve is equivalent to a spin structure on it.

Now, there is a geometrical characterization of SUSY-curves in
terms of superdivisors:
\proclaim{Theorem 4}   Let $\X$ be a supercurve of dimension
$(1,1)$. Then $\X$ is a SUSY-curve if and only if there is an
isomorphism of supercurves $\X\iso\X^c$ between $\X$ and the
supercurve of positive superdivisors of degree 1 (conjugate
fermions) on it inducing the identity on $X$. Moreover, there is a
one-to-one correspondence between such isomorphisms and spin
structures on $X$. \endproclaim
\proof
If $\O_{\X}=\O_X\oplus\Lc$, then the structure sheaf of
$\X^c$ is
$\O_X\oplus(\Lc^{-1}\otimes_{\O_X}\kappa_X)$, so that an
isomorphism $\X\iso\X^c$ inducing the identity on $X$
is nothing but a $\O_X$-module isomorphism
$\Lc^{-1}\otimes_{\O_X}\kappa_X\iso\Lc$, that is, an isomorphism
$\Lc\otimes_{\O_X}\Lc\iso\kappa_X$.\qed\enddemo

Theorem 3 and the former result, mean that for
SUSY-curves, $\Sc$-points are equivalent to relative positive
superdivisors of degree 1 on $\X\times\Sc\to\Sc$, as
Manin claimed in \cite{23}, and the
universal relative positive superdivisor of degree 1, gives in
this case nothing but Manin's superdiagonal:

Let $\X$ be a SUSY-curve.
If $\Delta$ denotes the ideal of the diagonal immersion
$\Delta\colon\X\hookrightarrow\X\times\X$, the kernel of the
composition $\Delta\to\Delta/\Delta^2\iso\Delta_\ast
\Omega^1_\X@>\partial>>\Delta_\ast Ber(\O_\X)$ is a
homogeneous ideal $\Cal I$ of $\O_{\X\times\X}$ thus defining a
sub-superscheme $\Delta^s$ called the {\sl superdiagonal\/}.

\proclaim{Lemma 3}  {\rm (Manin, \cite{23})} The superdiagonal
$\Delta^s=(X,\O_{\X\times\X}/\Cal I)$ is a closed sub-super\-scheme
of codimension $(1,0)$. In conformal coordinates $(z,\theta)$, it
can be described by the equation $$ z_1-z_2-\theta_1\theta_2=0 $$
where as usual $z_1=1\otimes z$ and $z_2=z\otimes
1$.\qed\endproclaim
According to Lemma 2, the superdiagonal is a
positive superdivisor. A simple local computation shows
that actually we have:
\proclaim{Theorem 5}   Let $\X$ be a SUSY-curve,
$\psi\colon\X\iso\X^c$ the natural isomorphism between
$X$ and the supercurve of positive superdivisors of degree 1 (conjugate
fermions,) and
$1\times\psi\colon\X\times\X\iso\X\times\X^c$ the induced
isomorphism. Then $$ \Delta^s=(1\times\psi)^{-1}(\Zc_1)\,,
$$
that is, the isomorphism $\psi\colon\X\iso\X^c$
given by the spin structure transforms by inverse image the
universal positive superdivisor of degree 1 into Manin's
superdiagonal.
\qed\endproclaim
\subheading{3. The superscheme of positive superdivisors of degree
$g$}

Let $\X$  be a smooth proper supercurve of
dimension $(1,1)$ as above.
\definition{Definition 9} The superscheme of positive superdivisors
of degree $g$ of $\X$ is the supersymmetric product $S
^g\X^c$ of the supercurve $\X^c$ of
positive superdivisors of degree 1.
\enddefinition
The universal superdivisor $\Zc^u_g$ of $\X\times
S^g\X^c$ is constructed as follows: let us
consider the natural projections
$$\align \pi_i\colon\X\times\X^c
\times\overset g\to\cdots\times\X^c
&\to\X\times\X^c\\
(x,x^c_1,\dots,x^c_g) &\mapsto (x,x^c_i)\,,\endalign$$
the positive superdivisors of degree 1, $\Zc_i=\pi
_i({\Zc}^u_1)\subset\X\times(\prod _{i=1}^g
{\X}^c)$ and the positive superdivisor of degree $g$, $ \Zc=
\Zc_1+\cdots +\Zc_g$.
\proclaim{Lemma 4}   There exists a unique positive relative
superdivisor ${\Zc}^u_g$ of degree $g$ of $\X\times S
^g\X^c\to S
^g\X^c$, such that $\pi^{\ast}({\Zc}^u_g)=\Zc$, where $\pi$ is the
natural morphism
 $$\pi\colon
\X\times(\prod _{i=1}^g {\X}^c)\to\X\times S^g
{\X}^c\,.
$$
\endproclaim
\demo{Proof} One has only to prove that  ${\Zc}^u_g=\pi
(\Zc)$ is the desired superdivisor. This can be done locally,
so that we can assume that $\X=\Spec A$ is affine and the line
bundles  $\Lc$ and $\kappa_X$ are trivially generated
respectively by $\theta$ and $dz$. Then, the local
equation of
${\Zc}^u_1$ is
$z\otimes 1-1\otimes
z-\theta\otimes\theta^c=0$ (see equation (4),)
and $\Zc$ is the superdivisor defined by the equation
$$
0=\prod _{i=1}^g( z-z_i-\theta{\theta}^c_i
)=z^g-(s_1+\theta\cdot\varsigma _1)z^{g-1}+\cdots
+(-1)^g(s_g+\theta\cdot\varsigma _g)\,,
$$
where $z_i=\pi
_i^{\ast}(1\otimes z),\theta^c_i =\pi
_i^{\ast}(1\otimes\theta^c)$ and $s_i,\varsigma _i$
are the even and odd symmetric functions corresponding to $z$
and $\theta^c$ (see Corollary 1.)
It follows that this last equation is also the local equation of
${\Zc}^u_g$ in $\X\times S^g{\X}^c$ and one
can readily check that $\pi
^{\ast}({\Zc}^u_g)=\Zc$.\qed\enddemo
\subheading{4. The representability theorem}

This paragraph will justify the above definitions by displaying
the representability theorem
\proclaim{Theorem 6}   The pair $(S^g{\X}^c,{\Zc}^u_g)$ represents
the functor of relative positive superdivisors of degree $g$ of
$\X$, that is, the natural map: $$
\align
\phi\colon\Hom(\Sc,S^g{\X}^c)&\to
\Div_\Sc^g(\X\times \Sc)\\ f&\mapsto (1\times
f)^{\ast}{\Zc}^u_g\,,
\endalign
$$
is a functorial isomorphism for every superscheme $\Sc$.
\endproclaim
\demo{Proof}

1) $\phi$ {\sl is injective:}

Let $
U=\Spec A\subset
 X$ be an open subscheme of the underlying ordinary curve $X$,
such that $\kappa_X$ and $\Lc$ are trivial generated
respectively by  $dz$, $\theta$. Let us consider the affine open
sub-superschemes $\U=\Spec (A\oplus \theta\cdot A)\hookrightarrow\X$
and ${\U}^c=\Spec (A\oplus\theta^c\cdot
A)\hookrightarrow{\X}^c$, where
$\theta^c=dz\otimes\omega_{\theta }\in\Gamma
(U,\kappa_X\otimes\Lc^{-1})=\Gamma (U,\Lc^c)$.

Now, $S^g{\U}^c\hookrightarrow S^g{\X}^c$ is an affine
open sub-superscheme and the symmetric functions $s_i(z)$,
$\varsigma _i(z,\theta^c)$ ($i=1,\dots ,g$) is a graded
system of parameters for the graded ring $S
^g_k(A\oplus\theta^c\cdot A)$. Let us denote it  simply by
$s_i,\varsigma_i$.

The family of the affine open sub-superschemes $S^g{\U}^c$ so
obtained (when $U$ ranges on the affine open subschemes of $X$
where $\kappa_X$ and $\Lc$ are trivial) is an open covering of
$S^g{\X}^c$ by affine open sub-superschemes such
that the universal positive superdivisor of
$\U\times S ^g {\U}^c\to\U$ is the closed sub-superscheme
${\Zc}^u_{\U}$ defined by the equation $$
z^g-(s_1+\theta\cdot\varsigma _1)z^{g-1}+\cdots
+(-1)^g(s_g+\theta\cdot\varsigma _g)=0\,.
$$

Then one has that for these affine open sub-superschemes the map
$$\align
\phi_{\U}\colon\Hom (\Sc,S^g {\U}^c)&\to
\Div_\Sc^g(\U\times\Sc)\\ f&\mapsto (1\times f)^{\ast}
{\Zc}^u_{\U}\,, \endalign
$$
is injective: In fact, we can assume that $\Sc$ is affine $\Sc=\Spec
B$. Now, the morphisms $f\colon \Sc\to S^g {\U}^c$ are
determined by the inverse images of the symmetric
functions $s_i$, $\varsigma _i$. But these inverse images are
determined by $(1\times f)^{\ast} {\Zc}^u_{\U}$ since
the coefficients of the characteristic polynomial of  $z\otimes
1$ acting by multiplication on the $B[\theta ]$-module
$\O_{(1\times f)^{\ast} {\Zc}^u_{\U}}$ are
$(-1)^i(f^{\ast}(s_i)+ \theta\cdot f^{\ast}(\varsigma _i))$.
This allows us to conclude.

A straightforward consequence of this fact is that the map
$\phi$ of the statement is injective for every superscheme
$\Sc$.

2) $\phi$ {\sl is an epimorphism:}

It is sufficient to prove that given a relative
positive superdivisor of degree $g$, $\Zc\subset\X\times
\Sc\to\Sc$, for every geometric point $p\in\Sc$ there
exist an open neighbourhood,
$\V\subset\Sc$, and a morphism $f_{\V}\colon\V\to S
^g {\X}^c$ such that $(1\times
f_{\V})^{\ast}({\Zc}^u_g)=\Zc\cap (\X\times
\V)=\Zc_{\V}$, for, in that case, these morphisms define a
morphism $f\colon\Sc\to S ^g
{\X}^c$ fulfilling $(1\times f)^{\ast}( {\Zc}^u_g)= \Zc$
by virtue of the former paragraph. Let $\pi\colon \X\times
\Sc\to\Sc$ be the natural projection and $
U=\Spec A\subset X$ an affine subscheme where $\kappa_X$
and $\Lc$ are trivial and such that (with the notation of the
beginning of this section) the affine open sub-superscheme
$\U\subset\X$  contains the superdivisor $\pi^{-1}(p)\cap
\Zc\hookrightarrow\X$. Then, $\W=\Sc-\pi(\Zc-\Zc\cap (\U\times
\Sc))$ is open, because $\pi$ is a proper morphism, and it
contains the point $p\in\Sc$. Let $\V=\Spec
B\subset\Sc$ be an affine open sub-superscheme containing $p$ and
contained in $\W$. By construction, if we put
$\Zc_{\V}=\Zc\cap\pi^{-1}(\V)$, then $\Zc_{\V}$
is a relative positive divisor of degree
$g$ of $\U\times \V\to\V$, so that it is affine $\Zc_{\V}=\Spec C$.
Let $dz$, $\theta$ be generators of $\kappa_X$ and
$\Lc$, as usual. Now, according to the definition of
superdivisor, the ring $C$ of $\Zc_{\V}$ is a
locally free module over $B[\theta]$ of rank $g$ and $\overline
C=C/\theta\cdot C$ is the ring of an ordinary divisor of degree
$g$ of $U\subset X$.

Let us consider the morphism
$f_{\V}\colon\Spec B=\V\to S^g{\U}^c=\Spec S
_k^g(A\oplus{\theta}^c\cdot A)$ induced by the ring
morphism $f^{\ast}_{\V}\colon S
_k^g(A\oplus{\theta}^c\cdot A)\to B$ defined, by
means of the determinant morphism, as follows: Let $S ^g_kA\to
B$ be the determinant morphism  defined by the quotient ring
$\overline C$ of $A\otimes_k B$. This morphism endows
$B$ with a structure of  $S ^g_kA$-algebra. But, by Lemma
1, one has $S^g_k(A\oplus {\theta}^c\cdot
A)=\medwedge_{S^g_kA} M$  for a certain free
$S^g_kA$-module
 $M$ generated by the odd symmetric functions $\varsigma _i$;
then, by the universal property of the exterior algebra,
defining $f^{\ast}_{\V}$  is equivalent to giving a
homogeneous morphism of degree zero of $S^g_kA$-modules, $M\to
B$. This morphism is actually characterized by the images of
the functions
$\varsigma _i$ ($i=1,\dots ,g$,) and we define these images
as the odd coefficients of the characteristic polynomial of
 $z\otimes 1$ acting on the $B[\theta ]$-module $C$ by
multiplication; this means that if the characteristic polynomial
is $z^g-(a_1+\theta\cdot b_1
)z^{g-1}+\cdots +(-1)^g(a_g+\theta\cdot b_g)$, then we define
$f^{\ast}_{\V}(\varsigma_i)=b_i$.

Moreover, one also has that $f^{\ast}_{\V}(s_i)=a_i$ and then
$(1\times f)^{\ast}_{\V}( {\Zc}^u_\U)$ is the relative
positive superdivisor of degree $g$ of $\U\times \V\to\V$
defined by the equation
$$z^g-(a_1+\theta\cdot b_1
)z^{g-1}+\cdots +(-1)^g(a_g+\theta\cdot b_g)\,.
$$
On the other hand, this is the characteristic polynomial
of $z\otimes 1$ acting by multiplication on the structure ring
of $\Zc_{\V}$, so that this polynomial vanishes on $ \Zc_{\V}$,
which means that $\Zc_{\V}$ is contained in $(1\times
f)^{\ast}_{\V}( {\Zc}^u_\U)$. Since both positive
superdivisors have the same degree, they are equal, thus
finishing the proof. \qed\enddemo
\subheading{5. The case of SUSY-curves}

If $\X$ is a SUSY-curve, there exists an
isomorphism $\psi\colon\X\iso\X^c$ between $X$ and the
supercurve of positive superdivisors of degree 1, as we proved in
subsection 4.1. Then we have an
isomorphism $S^g\X\to S^g\X^c$ between the supersymmetric
product of $\X$ and the superscheme $S^g\X^c$ of positive
superdivisors of degree $g$ on $X$, so that the
representability theorem now reads (see \cite{8}):
\proclaim{Theorem 7}   Let $\X$ be a SUSY-curve. The supersymmetric
product
 $S^g\X$ represents the functor of positive superdivisors on
$\X$, that is, there exists a universal relative positive
superdivisor $\Zc^u_g$ of degree $g$ of $\X\times S^g\X\to
S^g\X$ such that the natural map
$$
\align
\phi\colon\Hom(\Sc,S^g\X)&\to
\Div_\Sc^g(\X\times \Sc)\\ f&\mapsto (1\times
f)^{\ast}\Zc^u_g\,,
\endalign
$$
is a functorial isomorphism for every superscheme $\Sc$.
\endproclaim
Moreover, since
$1\times\psi\colon\X\times\X\iso\X\times\X^c$
transforms by inverse image the universal positive superdivisor
of degree 1 into Manin's superdiagonal, the universal
superdivisor of $\X\times S^g\X\to
S^g\X$ for SUSY-curves is constructed as in Lemma
4 with Manin's superdiagonal playing the role of
${\Zc}^u_1$.

Summing up, only for SUSY-curves, `unordered families of $g$
points' (the points of $S^g\X$) are
equivalent to `superdivisors of degree $g$' (the points of
$S^g\X^c$.) \medskip
\noindent
{\bf Acknowledgments.} We thank J.M. Mu\~noz Porras for many
enlightening comments about Jacobian theory and the geometry of
 the symmetric products, J.M. Rabin for drawing references
\cite{28} and \cite{29} to our attention and J. Mateos
Guilarte who first introduced us to the vortex equations. We also
thank the anonymous referee for some helpful suggestions to
improve the original manuscript.

\enddocument \bye